\input amstex
\documentstyle{amsppt}
\magnification1200
\NoBlackBoxes
\pagewidth{6.5 true in}
\pageheight{9.25 true in}
\document

\topmatter
\title 
Two $S$-unit equations with many solutions 
\endtitle
\author S. Konyagin and K. Soundararajan 
\endauthor 
\thanks The second author is partially supported by the National Science 
Foundation.
\endthanks
\address 
{Department of Mathematics and Mechanics, Moscow State University, 
Moscow 119992, Russia} 
\endaddress
\email{konyagin{\@}ok.ru}\endemail
\address{Department of Mathematics, University of Michigan, Ann Arbor, 
Michigan 48109, USA} \endaddress
\email{ksound{\@}umich.edu} \endemail
\endtopmatter

\document

\head 1. Introduction \endhead 

\noindent In this note we consider two $S$-unit equations 
for which we will exhibit many solutions.  
Our first problem concerns solutions to the equation $a+b=c$ where 
$a$, $b$, and $c$ are coprime integers such that 
all prime factors of $abc$ lie in a given set $S$ of $s$ primes.  
In [8] J.-H. Evertse showed that this 
$S$-unit equation has at most $\exp(4s+6)$ solutions.  On the other 
hand, in [6] P. Erd{\H o}s, C. Stewart, and R. Tijdeman 
showed that there exist arbitrarily large sets $S$ such that the $S$-unit equation 
$a+b=c$ has more than $\exp((4-\epsilon) \sqrt{s/\log s})$ solutions (see 
also [9] for a refinement of their result).  
The set $S$ that they exhibited is rather special, and they conjectured 
that if $S$ were the set of the first $s$ prime numbers then there 
should be  $\gg \exp(s^{\frac 23 -\epsilon})$ solutions 
to the $S$-unit equation.  Moreover, for any set $S$ they conjectured 
that there are $\ll \exp(s^{\frac 23+\epsilon})$ solutions.  We 
remark that recently J. Lagarias and K. Soundararajan [12] have shown that 
if $S$ is the set of the first $s$ prime numbers and the Generalized 
Riemann Hypothesis is true then the $S$-unit equation has 
$\gg \exp(s^{\frac 1{12}-\epsilon})$ solutions.  Our first result 
improves the construction of Erd{\H o}s, Stewart, and Tijdeman 
and shows the existence of arbitrarily large sets $S$ with 
more than $\exp(s^{2-\sqrt{2}-\epsilon})$ solutions.

\proclaim{Theorem 1}  Let $\beta$ be any positive 
number with $\beta < 2-\sqrt{2}$.  There exist arbitrarily 
large sets $S$ of $s$ prime numbers such that the 
$S$-unit equation $a+b=c$ has at least $\exp(s^{\beta})$ 
solutions in coprime integers $a$, $b$ and 
$c$ having all their prime factors from $S$.
\endproclaim 

The second $S$-unit equation that we will consider is a special 
case of the first: namely, the equation $a+1=c$ with all prime 
factors of $ac$ lying in the set $S$.  Although this is a much more 
restrictive equation than our first, we are able to find arbitrarily 
large sets $S$ with many solutions to this equation.  

\proclaim{Theorem 2}  There exist arbitrarily large sets $S$ 
of $s$ prime numbers such that the equation $a+1=c$ has 
at least $\exp(s^{\frac 1{16}})$ solutions where all prime factors of $ac$
lie in $S$.  In fact, there exist arbitrarily large integers 
$N$ such that 
$$
\# \{ d: \ \ d(d+1) | N\} \ge \exp((\log N)^{\frac 1{16}}).
$$
\endproclaim 

The second, stronger, conclusion of Theorem 2 advances a line 
of inquiry initiated by Erd{\H o}s and R.R. Hall [5].  They showed the 
existence of arbitrarily large numbers $N$ with $\#\{ d: \ \ d(d+1) |N\} 
\gg (\log N)^{\sqrt{e}-\epsilon}$.  From the work of A. Hildebrand [10] 
on consecutive smooth numbers it follows that there are large 
$N$ with $\#\{ d: \ \ d(d+1) |N\} 
\gg (\log N)^A$ for any given positive number $A$. In [1] A. Balog, 
Erd{\H o}s, and G. Tenenbaum quantified this and obtained large $N$ with 
$\#\{ d: \ \ d(d+1) | N\} \gg (\log N)^{\log_3 N/9\log_4 N}$ where 
$\log_3$ and $\log_4$ denote the third and fourth iterated logarithms.   For upper bounds 
on the quantity $\# \{d: \ \ d(d+1)|N\}$ we refer the reader to [3], [4], and [7]. 
  
There are at least $x^{\frac {1+\delta}{2+\delta}+o(1)}=x^{\frac 12+\frac{1}{4+2\delta}+o(1)}$ 
square-free numbers below $x$ all of whose prime factors lie below $(\log x)^{2+\delta}$.  
If these numbers were randomly distributed then we would expect to find 
about $x^{\frac{1}{2+\delta}+o(1)}$ pairs of such consecutive numbers.  This 
suggests that there should be arbitrarily large $N$ with $\# \{ d: \ \ d(d+1)|N\} 
\ge \exp((\log N)^{\frac 12-\epsilon})$.  We venture the guess that 
for any set $S$, the $S$-unit equation $a+1=c$ has no more than 
$\exp(s^{\frac 12+\epsilon})$ solutions, but nothing substantially better than Evertse's 
bound appears to be known.
 
 This work was completed when both authors were visiting the 
 Centre de Rechereches Math{\' e}matiques (CRM) in Montr{\' e}al.  
 We are grateful to the CRM for their support and excellent working conditions.  We 
 are also grateful to Antal Balog and Andrew Granville for their interest and 
 encouragement.  
  
\head 2.  Proof of Theorem 1 \endhead 

\noindent Let $y$ be a large real number and let 
$\beta$ and $\gamma$ be  real numbers in $(0,1)$. Consider the 
set ${\Cal L}$ which consists of square-free numbers 
$\ell$ having exactly $[y^{\beta}]$ prime factors 
each from the interval $[y/2,y]$.  Consider also the 
set ${\Cal M}$ which contains square-free numbers $m$ 
having exactly $[\gamma y^{\beta}]$ prime factors each 
from the interval $[y/4,y/2)$.  Note that 
the elements of ${\Cal L}$ are coprime to elements of ${\Cal M}$.  Further note
that 
$$
|{\Cal L}| = \binom{\pi(y)-\pi(y/2)}{[y^{\beta}]} 
= L^{1-\beta+o(1)}, 
$$
where $L = y^{[y^{\beta}]}$, and similarly 
$$
|{\Cal M}| = L^{\gamma(1-\beta) +o(1)}. 
$$

Pick a number $m\in {\Cal M}$ and let $r({\Cal L};a,m)$ 
denote the number of elements of ${\Cal L}$ lying 
in the residue class $a \pmod m$.  By Cauchy-Schwarz 
we know that 
$$
\sum_{a=1}^{m} r({\Cal L};a,m)^2 \ge \frac 1m \Big(\sum_{a=1}^{m} 
r({\Cal L};a,m)\Big)^2  = \frac{|{\Cal L}|^2}{m}.
$$
The left hand side counts the pairs $(\ell_1,\ell_2)$ with $\ell_1
\equiv \ell_2 \pmod m$.  This congruence has $|{\Cal L}|$ trivial 
solutions, and if $m <|{\Cal L}|/2$ then we are 
guaranteed $\gg |{\Cal L}|^2/m$ non-trivial solutions.  
Since each element of ${\Cal M}$ is below $y^{[\gamma y^\beta]} 
\le y L^{\gamma}$ we conclude that if $\gamma < 1-\beta$ 
then there exist $\gg L^{2(1-\beta)-\gamma+o(1)}$ 
non-trivial pairs $(\ell_1,\ell_2)$ with $\ell_1\equiv \ell_2 
\pmod m$.  Therefore, if $\gamma <1-\beta$ there exist 
$\gg L^{2(1-\beta)-\beta\gamma +o(1)}$ triples $(m,\ell_1,\ell_2)$ 
with $m\in {\Cal M}$, $\ell_1 \neq \ell_2\in {\Cal L}$ and 
$\ell_1\equiv \ell_2\pmod m$.

Suppose below that $\gamma <1-\beta$ and consider the 
ratios $(\ell_1-\ell_2)/m$ arising from the triples produced 
above.  Restricting to positive ratios, we have produced 
$\gg L^{2(1-\beta)-\beta\gamma+o(1)}$ such ratios, all 
below $L^{1-\gamma +o(1)}$.  Therefore if 
$2(1-\beta)-\beta\gamma > 1-\gamma$ then we can find a popular 
number $u\le L^{1-\gamma +o(1)}$ which occurs 
as a ratio more than $L^{2(1-\beta)-\beta\gamma +\gamma-1+o(1)}$ 
times.  

Summarizing, we see that if $\gamma < 1-\beta$ and $(2+\gamma)(1-\beta)>1$ 
then there is a number $u\le L^{1-\gamma+o(1)}$ such that the 
equation $\ell_1= \ell_2 + mu$ has more than $L^{(2+\gamma)(1-\beta)-1+o(1)}$ 
solutions in integers $\ell_1 \neq \ell_2 \in {\Cal L}$ and $m\in {\Cal M}$.  
We already know that $\ell_1$ and $\ell_2$ are coprime to $m$, so if $\ell_1$ and
$\ell_2$ have a common factor then it must be a divisor of $u$.  Since there 
are at most $L^{o(1)}$ divisors of $u$, after removing common factors, we find 
that for some divisor $v$ of $u$, the equation $\ell_1=\ell_2 +v m$ has 
$\gg L^{(2+\gamma)(1-\beta)-1+o(1)}$ solutions in coprime integers $\ell_1$, $\ell_2 \in{\Cal L}$, 
and $m\in {\Cal M}$.  Take $S$ to be the set of all primes in $[y/4,y]$ union 
the prime factors of $v$.  Then $|S| \le \pi(y)-\pi(y/4) + \log v \le y$, 
and we have exhibited more than $\exp(y^{\beta})$ solutions to this $S$-unit 
equation.  If $\beta < 2-\sqrt{2}$ then we can find a $\gamma$ satisfying the 
conditions $\gamma <1-\beta$ and $(2+\gamma)(1-\beta)>1$, and so Theorem 1 follows.

\head 3.  Proof of Theorem 2 \endhead 

\noindent Throughout we let $y$ be a large real number.  
We need first the following zero-density result which may 
be found in [11] (see the Grand Density Theorem 10.4 on page 260). 

\proclaim{Lemma 3.1} There exists a constant $C>0$ such that 
for any $\frac 12 \le \alpha < 1$ the region 
$$
{\Cal R}(\alpha,y):= \{ s: \ \ \text{Re }(s) \ge \alpha , \ \ 
|\text{Im }(s)|\le y\},
$$
contains at most $(Q^2 y)^{C(1-\alpha)+o(1)}$ zeros 
of primitive Dirichlet $L$-functions with conductor below $Q$.
It is permissible to take $C=\frac {12}{5}$.   
\endproclaim

\proclaim{Proposition 3.2}  Let $\beta$ be a real number with 
$0< \beta < 1-3C(1-\alpha)$.  Let $K=[y^{\beta}]$ and put $Z =y^K$.
There exist $\gg Z^{1-\beta +o(1)}$ square-free numbers $q$ having 
exactly $K$ prime factors each from the interval $[y/2,y]$, 
and such that for every non-trivial character $\pmod q$ 
the corresponding $L$-function has no zeros in ${\Cal R}(\alpha,y)$. 
\endproclaim 

\demo{Proof} Clearly there there are $\binom{\pi(y)-\pi(y/2)}{K}$ 
square-free integers $q$ having exactly $K$ prime 
factors each from the interval $[y/2,y]$.  We must exclude 
those moduli for which there exists a non-trivial character 
whose $L$-function has a zero in ${\Cal R}(\alpha,y)$.  A bad modulus 
$q$ must be divisible by some number $d$ with $j$ prime factors (so $(y/2)^j \le 
d\le y^j$ and $1\le j\le K$) such that there is a primitive character $\mod d$ 
whose $L$-function has a zero in ${\Cal R}(\alpha,y)$.  By 
Lemma 3.1 there are at most $y^{(2j+1)(C(1-\alpha)+o(1))}$ 
possibilities for $d$.  Given a $d$ there are at most 
$\binom{\pi(y)-\pi(y/2)}{K-j}$ multiples of $d$ that 
must be excluded.  Thus we must exclude at most 
$$
\sum_{j=1}^K y^{(2j+1)(C(1-\alpha)+\epsilon)} \binom{\pi(y)-\pi(y/2)}{K-j}
$$
moduli.  Since $\beta< 1-3(1-\alpha)$ this is small compared to $\binom{\pi(y)-\pi(y/2)}{K}$ 
and so we have $\gg \binom{\pi(y)-\pi(y/2)}{K} = Z^{1-\beta+o(1)}$ suitable 
moduli $q$.

\enddemo 

\proclaim{Proposition 3.3} Let $X=Z^{\gamma}$ and 
suppose that $\gamma (1-\alpha-\beta) > 1$.  Let $q$ 
be one of the moduli produced in Proposition 3.2.  Then 
there are $\gg Z^{(1-\beta)\gamma -1+o(1)}$ integers 
$\ell \le X$ with each $\ell$ being square-free, 
divisible only by primes below $y$, and $\ell \equiv 1\pmod q$. 
\endproclaim 

Assuming this Proposition for the moment we show how to 
deduce Theorem 2.

\demo{Proof of Theorem 2} Let $\alpha$, $\beta$, and $\gamma$ 
be as in Lemma 3.1, Propositions 3.2 and 3.3.  That is 
$$
\tfrac 12\le \alpha < 1, \qquad 
0< \beta < 1-3C(1-\alpha), \qquad  
\text{and } \gamma (1-\alpha-\beta) >1. 
\tag{3.1}
$$  
By Propositions 3.2 and 3.3 
we know that there are at least $Z^{(1-\beta)(1+\gamma)-1+o(1)}$ 
pairs $(\ell,q)$ satisfying the conclusions of those Propositions.  
Consider the ratio $(\ell-1)/q$ which is an integer which lies 
below $2^K X/Z < Z^{\gamma -1 +o(1)}$.   If 
$$
(1-\beta)(1+\gamma)-1 > \gamma -1, \tag{3.2}
$$ 
then there is a popular value $m$ 
which occurs as the ratio $(\ell-1)/q$ at least $Z^{1-\beta-\beta \gamma +o(1)}$ 
times.  Take $N=m \prod_{p\le y} p$ and note that 
if $(\ell-1)/q=m$ then $qm$ and $\ell =qm+1$ are consecutive divisors of $N$.  Therefore 
$$
\#\{d: \ \ d(d+1) | N\} \ge Z^{1-\beta -\beta \gamma +o(1)} 
\ge \exp((\log N)^{\beta}),
$$
since by the prime number theorem $N=e^{y+o(y)}$, and $\log Z 
= (1+o(1))y^{\beta} \log y$.

To complete the proof we need only find the largest $\beta$ for which (3.1) and (3.2) 
hold.  A little calculation shows that it is best to take 
$\gamma$ slightly larger than $3C+\sqrt{9C^2 +3C}$, take $\alpha= 1- \frac{1+1/\gamma}{3C+1}$, 
and $\beta$ is then slightly smaller than $(1+3C + \sqrt{9C^2 + 3C})^{-1}$.  Since 
$C=\frac {12}{5}$ is permissible we conclude that $\beta=\frac{1}{16}$ is allowed. 

\enddemo 

It remains finally to prove Proposition 3.3.  To this end we require the 
following Lemma.  
\proclaim{Lemma 3.4}  Let $q \le Z$ be one of 
the moduli produced in Proposition 3.2 so that 
$L(s,\chi)$ has no zeros in the region ${\Cal R}(\alpha,y)$, and 
suppose that $\beta <1-\alpha$.  
For any complex number $s$ with Re$(s)>0$ we define 
$$
F(s,\chi;y) =\sum\Sb \ell=1\\ p|\ell \implies p\le y\endSb^{\infty}
\frac{\mu(\ell)^2 \chi(\ell)}{\ell^s} = \prod_{p\le y} \Big( 1+\frac{\chi(p)}{p^s}\Big).
$$
For any $\epsilon>0$, if $|t| \le y/2$ then we have 
$$
|F(\alpha +\epsilon+it,\chi;y)| \ll_{\epsilon} (qy)^{\epsilon}, 
$$
while if $|t|>y/2$ we have 
$$
|F(\alpha+\epsilon+it,\chi;y)| \ll \exp(y^{1-\alpha}).
$$
\endproclaim 

\demo{Proof} Taking logarithms it suffices to estimate $\sum_{p\le y} 
\chi(p)p^{-\alpha-\epsilon-it}$.  Since $\alpha <1$ this is 
trivially $\le y^{1-\alpha}$ and the second assertion follows.  

If $z\le y$ then note that 
$$
\align
\sum_{n\le z} \Lambda(n) \chi(n)n^{-it} 
&= \frac{1}{2\pi i} \int_{1+\frac{1}{\log y}-i\infty}^{1+\frac{1}{\log y}+i\infty} 
-\frac{L^{\prime}}{L} (w+it,\chi)\frac{z^{w}}{w} dw  
\\
&= -\sum\Sb \rho \\ |\rho -it| \le z/2 
\endSb \frac{z^{\rho-it}}{\rho-it} +O(\log^2 qy),
\\
\endalign
$$
by following closely the standard argument in prime number theory leading 
to the `explicit formula' for primes (see for example H. Davenport [2]); here 
$\rho$ runs over non-trivial zeros of $L(s,\chi)$.  By assumption Re$(\rho)\le \alpha$ 
for each zero counted in our sum.  Since there are $\ll \log qy$ zeros in each interval 
$k \le |\rho -it| \le k+1$ for $0\le k\le y$ we conclude that 
$$
\sum_{n\le z} \Lambda(n) \chi(n)n^{-it} \ll z^{\alpha} \log (qy) \log z + 
\log^2 (qy) \ll z^{\alpha} \log (qy) \log z.
$$
Trivially we also have that this sum is bounded by $\ll z$.  Using these two 
estimates and partial summation we easily deduce that 
$$
\sum_{2\le n\le z} \frac{\Lambda(n)\chi(n)}{n^{\alpha+\epsilon+it} \log n} 
\ll (\log qy)^{1-\frac{\epsilon}{1-\alpha}}.
$$
This proves the Lemma.

\enddemo 

\demo{Proof of Proposition 3.3} Using 
the orthogonality of characters $\pmod q$ we see that 
$$
\sum\Sb \ell \equiv 1 \pmod q \\ p|\ell \implies p\le y \endSb 
\mu(\ell)^2 e^{-\ell/x} 
= \frac{1}{\phi(q)}\sum\Sb (\ell, q)=1\\ p| \ell 
\implies p \le y \endSb e^{-\ell/x} + \frac{1}{\phi(q)}
\sum\Sb \chi \pmod q  \\ \chi \neq \chi_0 \endSb  
\sum\Sb p|\ell \implies p\le y \endSb \chi(\ell) \mu(\ell)^2 e^{-\ell/x}. 
\tag{3.3} 
$$

We now obtain an upper bound for the contribution from 
non-trivial characters to (3.3).  
For any $c>0$ we 
have 
$$
\sum\Sb p|\ell \implies p\le y\endSb \chi(\ell) \mu(\ell)^2 e^{-\ell/x} 
= \frac{1}{2\pi i} \int_{c-i\infty}^{c+i\infty} F(s,\chi;y) x^s \Gamma(s) ds.
$$
We take $c=\alpha+\epsilon$ and estimate the integral using Lemma 3.4.  
Since $|\Gamma(c+it)|$ decays exponentially in $|t|$ by Stirling's 
formula, we obtain that the above is 
$\ll x^{\alpha+\epsilon} (qy)^{\epsilon}$.  
Thus we conclude that 
$$
\sum\Sb \ell \equiv 1 \pmod q \\ p|\ell \implies p\le y \endSb 
\mu(\ell)^2 e^{-\ell/x} 
= \frac{1}{\phi(q)}\sum\Sb (\ell, q)=1\\ p| \ell 
\implies p \le y \endSb e^{-\ell/x} + O(x^{\alpha+\epsilon} (qy)^{\epsilon}). \tag{3.4}
$$

We take $x=X/\log X$ in (3.4) and note that 
$$
\sum\Sb \ell \le X \\  \ell \equiv 1 \pmod q \\ p|\ell \implies p\le y \endSb 
\mu(\ell)^2 
\ge \sum\Sb \ell \equiv 1 \pmod q \\ p|\ell \implies p\le y \endSb 
\mu(\ell)^2 e^{-\ell/x} + O(1).
$$
Now 
$$
\sum\Sb (\ell,q)=1\\ p|\ell \implies p\le y\endSb \mu(\ell)^2 e^{-\ell/x} 
\gg \sum\Sb \ell \le x \\ (\ell, q)=1\\ p|\ell \implies p\le y\endSb 
\mu(\ell)^2 \ge \binom{\pi(y)-\omega(q)}{[\log x/\log y]} 
= Z^{\gamma(1-\beta) +o(1)}.
$$
Using (3.4), and recalling that $q\le Z$ and the hypothesis 
that $\gamma(1-\beta)-1 > \gamma \alpha$, we obtain (choosing $\epsilon$ small enough) the Proposition.
\enddemo

\enddocument